\documentclass[preprint,12pt]{elsarticle}
\usepackage{amsmath,latexsym,amssymb, enumerate, amsthm, epsfig }
\usepackage{amsmath,amsfonts,amsthm,amscd,amsxtra,enumerate}
\usepackage[paperheight=297mm, paperwidth=210mm, left=30mm, right=20mm]{geometry}
\usepackage{hyperref}
\usepackage[hyperref]{xcolor}
\hypersetup{colorlinks=true,}

\begin{document}

\begin{frontmatter}

\title{Analytical and numerical solution of a transport equation for resonantly interacting waves in MHD for a van dar Waals gas}


\author{Harsh V. Mahara and V. D. Sharma}

 \address{Department of Mathematics, Indian Institute of Technology Bombay,
	Powai, Mumbai-400076}

\begin{abstract}
In this paper, we characterized an analytical and numerical study of the resonant interaction between waves in MHD. A system of evolution equations is derived; we focus on the study of the interaction between a selected triad. The resulting evolution equation contains a dispersive term in addition to the nonlinear term and convolution term. Effects of the influence of van der Waals parameter and magnetic field on the formation and structure of solitons are studied.

\end{abstract}


\end{frontmatter}


\section{Introduction}
The study of resonant interaction in weakly nonlinear waves has received increasing attention from mathematicians in the recent past \cite{MR338576, MR867874, MR1087089, MR1881084, MR1297667}. There has been widespread interest in the nonlinear phenomena mainly due to the so called evolution equations,
derived from a system of PDEs, representing an essential aspect of the original system \cite{ MR1166185, MR0426722, MR2259844, MR3623396, MR1013433}. To study the wave interactions in one dimension, Majda and Rosales \cite{MR760229} have derived a system of integro-differential equations and have shown their physical applications to the gas dynamics; the analytical and numerical solutions of this equation were studied by Majda, Rosales, and Schnobek \cite{MR975485}. The theory of resonant interaction was applied to elasticity and dispersive plasma in \cite{MR1160152} and \cite{MR0029616}, respectively. Ali and Hunter \cite{MR1616005} applied the Majda-Rosales theory to the MHD system, including the viscous and dispersive effects in one dimension with ideal gas background.  The MHD wave interaction shows dispersive behavior that is different from acoustic and elastic waves; this leads to the KdV-type term in the interaction equation. Zabusky and Kruskal \cite{MR716190} studied the solitary wave solutions of the KdV equation with weak dispersion.

In this chapter, we study the resonant interaction of weakly nonlinear waves in the MHD system with a van der Waals equation of state and derive an evolution equation corresponding to the fast magnetosonic entropy wave triad. The far-field behaviour of the underlying equation is studied analytically and numerically taking into account the effects of magnetic field and the van der Waals gas. The evolution equation has a Burgers equation type nonlinear term, a weak dispersive term, and a 
(weakly dispersive \cite{MR1719749}) convolution term corresponding to the interaction between waves. We notice the presence of solitons and study the effect of magnetic field  and the real gas parameter $ b $ on the evolution, shape and behavior of solitary wave profiles.

The work is organized as follows: The basic equations and formulation of the problem are given in Section  \ref{basic equation1}. The detailed derivation of the system of transport equations for the wave amplitudes, exhibiting nonlinearity, dissipation, and dispersion is given in Section  \ref{evo1}. The evolution equations related to the fast magnetosonic and entropy wave triad with certain assumption is developed in the Section
\ref{evo}. The numerical results, exhibiting the magnetic and real gas effects, are displayed in Section \ref{num1}. Finally, we concluded this chapter with a discussion of our results in Section \ref{con1}.

\section{Basic equations}\label{basic equation1}
For one dimensional motion, the equations of MHD can be written as a system in the  following form \cite{cab70}:
\begin{equation}\label{b.2.1}
\setlength{\jot}{5pt}
\begin{aligned}
&\rho_{t}+u\rho_{x}+\rho u_{x}=0,\qquad\qquad\qquad\qquad\qquad\qquad\qquad\qquad\qquad\qquad\qquad\qquad\qquad\\
&\rho(u_{t}+uu_{x})+ (P+\frac{1}{2}B^{2})_{x}=\frac{4}{3}\mu u_{xx},\\
&\rho(v_{t}+uv_{x})- B_{1}B_{2x}=\mu v_{xx},\\
&\rho(w_{t}+uw_{x})- B_{1}B_{3x}=\mu w_{xx},\\
&\rho(s_{t}+us_{x})= \kappa T_{xx} +\frac{4}{3}\mu u_x^2 + \mu (v_x^2 +w_x^2) + \eta (B_{2x}^2 + B_{3x}^2),\\ 
&B_{2t}+u B_{2x} + B_{2}u_{x}-B_{1}v_{x}=\eta B_{2xx} + \chi B_{1}B_{3xx},\\
&B_{3t}+u B_{3x} + B_{3}u_{x}-B_{2}w_{x}=\eta B_{3xx} - \chi B_{1}B_{2xx},\\
&B_{1t}=B_{1x}=0,
\end{aligned}
\end{equation}
where, $\, P $ is the pressure, $ \rho $ is the fluid density, $ s $ is the entropy, $ (u,v,w) $ is the fluid velocity vector, $ \mu $ is the viscosity,
$ \kappa $ the thermal conductivity, $ T $ the temperature, $ \eta $ the magnetic diffusion, $ \chi $ the hall parameter, 
$ (B_{1},B_{2},B_{3}) $  the magnetic field vector, and $ B^2 =B_{1}^2+B_{2}^2+B_{3}^3.$
The last equation of \eqref{b.2.1} implies that $ B_{1}(x,t)= \text{constant}, $ which reduces the number of the equations to seven, and is due to the fact the magnetic field is divergence free. We have considered the direction of propagation along the $ x $ axis, and the subscripts $ x $ and $ t $ denote the partial differentiation
with respect to the respective variable.

The system of  equations \eqref{b.2.1} is supplemented by an equation of state which in our case is the van der Waals equation of state of the form \cite{MR1419777}

\begin{equation}\label{b.2.2}
P = K_{0}\delta\frac{\rho^{1+\delta}\exp(\delta s/R)}{(1-b\rho)^{1+\delta}},
\qquad
T = K_{0}\delta\frac{\rho^{\delta}\exp(\delta s/R)}{R(1- b\rho)^{\delta}},
\end{equation}
where, $ K_{0} $ is a constant, $ \delta $ is a dimensionless material dependent quantity defined as $ \delta = R/c_{v} $ with $ c_{v} $ the specific heat at constant volume and $ R $ the specific gas constant;
$ \delta $ lies in the interval $ 0 < \delta \leq 2/3  $ with $  \delta= 2/3  $, for a monoatomic fluid and the parameter $ b $ represents the van der Waal excluded volume.

The system \eqref{b.2.1} can be written in the vector matrix notation as
\begin{equation}\label{b.2.3}
\mathbf{U}_{t}+\mathbf{A}(\mathbf{U})\mathbf{U}_{x}= \mathbf{M}(\mathbf{U})\mathbf{U}_{xx}+
\mathbf{N}(\mathbf{U})[\mathbf{Q} (\mathbf{U})]_{x}\mathbf{U}_{x},
\end{equation}
where $ \mathbf{U}=(\rho,u,v,w,s,B_{2},B_{3})^{tr} $ is the state column vector; $ \mathbf{A},\,\mathbf{M},\,\mathbf{N},\;\text{and}\;\mathbf{Q} $ are square matrices
of order 7 having components $ A_{ij},\,N_{ij},\,M_{ij},\;\text{and} \;Q_{ij} $ ,respectively; the non-zero components are as follows:
\begin{equation*}\label{b.2.4}
\setlength{\jot}{5pt}
\begin{aligned}
& \quad\qquad 
A_{ij}=u,~ i=1,\ldots,7, \quad~
A_{12}=\rho,\quad~ 
A_{21}=\frac{P_{\rho}}{\rho},\quad~ 
A_{25}=\frac{P_{s}}{\rho},\quad~ 
A_{26}=\frac{B_{1}}{\rho},
\\
& \quad\qquad 
A_{27}=\frac{B_{2}}{\rho},\quad~ 
A_{36}=\frac{-B_{1}}{\rho},\quad~ 
A_{47}=\frac{-B_{1}}{\rho},\quad~ 
A_{62}=B_{2},\quad~ 
A_{63}=-B_{1},\quad~ 
\\
& \quad\qquad 
A_{72}=B_{3},\quad~ 
A_{74}=-B_{1},\quad~ 
M_{22}=\frac{\mu'}{\rho},\quad~ 
M_{33}=\frac{\mu}{\rho},\quad~ 
M_{44}=\frac{\mu}{\rho},\quad~ 
\\
& \quad\qquad 
M_{51}=\frac{\kappa T_{\rho}}{\rho T},\quad~ 
M_{51}=\frac{\kappa T_{s}}{\rho T},\quad~ 
M_{66}=\eta,\quad~ 
M_{67}=-\chi B_{1},\quad~ 
M_{77}=\eta,\quad~ 
\\
& \quad\qquad 
N_{51}=\frac{\kappa}{\rho T},\quad~ 
N_{52}=\frac{\mu'}{\rho T},\quad~ 
N_{53}=\frac{\mu}{\rho T},\quad~ 
N_{56}=\frac{\mu}{\rho T},\quad~
Q_{11}=T_{\rho},\quad~ 
\\
& \quad\qquad 
Q_{15}=T_{s},\quad~
Q_{22}=u,\quad~
Q_{33}=u,\quad~ 
Q_{34}=w,\quad~
Q_{66}=B_{2},\quad~
Q_{67}=B_{3}.\quad~
\\
\end{aligned}
\end{equation*}

In the system \eqref{b.2.3}, matrices $ \mathbf{M} ,$  $ \mathbf{N} ~\text{and}~ \mathbf{Q} $ correspond to the dispersive and diffusive parts; we can split $\mathbf{M}$ as $ \mathbf{M} = \mathbf{M}_{d} + \mathbf{M}_{v}, $ where the dispersive 
part $ \mathbf{M}_{d} $ is proportional to $ \chi $ and the diffusive part $ \mathbf{M}_{v} $ is proportional to 
$\mu,\,\kappa,~\text{or}~\eta  .$

\section{Derivation of evolution equations}\label{evo1}
It is to be noted that the system \eqref{b.2.1} is parabolic in nature with the second-order spatial derivatives. With the removal of viscosity, thermal conduction, magnetic diffusion, and Hall effect, the system \eqref{b.2.1} reduces to a hyperbolic system. To study the interaction between the MHD waves with real gas effects, we derive a system of evolution equations in the present section.

We analyse the interaction of waves, which propagate through a constant background state
$ \mathbf{U}_{0} = (\,\rho_{0},\, 0,\, 0,\, 0,\, s_{0},\, B_{02},\, 0\, )^{tr} .$
Since the left-hand side of the system \eqref{b.2.3} is  hyperbolic, it admits seven families of characteristic velocities at $  \mathbf{U}=\mathbf{U}_{0} $ given by
\begin{equation}\label{b.2.5}
\setlength{\jot}{5pt}
\begin{aligned}
&\lambda_{1}=0,\qquad\qquad\qquad\qquad\qquad\qquad\qquad\qquad\qquad\qquad\qquad\qquad\qquad\\
&\lambda_{2,3}=\mp \frac{B_{01}}{\sqrt{\rho_{0}}},\\
&\lambda_{4,5}=\mp \sqrt{\frac{1}{2}\left(c_{0}^2 +
	\frac{1}{\rho_{0}}(B_{01}^2+B_{02}^2)\right) 
	-\sqrt{c_{0}^2+\frac{1}{\rho_{0}}(B_{01}^2+B_{02}^2-4c_{0}^2B_{01}^2)}},\\
&\lambda_{6,7}=\mp \sqrt{\frac{1}{2}\left(c_{0}^2 +
	\frac{1}{\rho_{0}}(B_{01}^2+B_{02}^2)\right) +\sqrt{c_{0}^2+\frac{1}{\rho_{0}}(B_{01}^2+B_{02}^2-4c_{0}^2B_{01}^2)}},\\
\end{aligned}
\end{equation}
where $ c_{0}=( P_{\rho_{0}})^{1/2} $ is the sound speed. The wave with the speed $ \lambda_{1} $ corresponds to the convection of the entropy with the particle velocity and is called entropy wave, while waves with speeds $ \lambda_{2,3},\, \lambda_{4,5},\, \lambda_{6,7}\, $ correspond to the left and right moving Alfvén waves, slow magnetoacoustic and fast magnetoacoustic waves and are denoted by $ c_{a},c_{s},~\text{and}~c_{f} $, respectively.
In certain degenerate cases some wave  speeds coincide; for instance, when 
$ B_{01}=0,\;B_{02}\neq 0, $ we have $ c_s=c_a=0 $ as an eigenvalue of multiplicity five;
when $ B_{01}\neq0,\; B_{02}= 0,\;\text{and}\; c_{a}=c_{s} $ we have $ c_{s}=c_{a}=c_{f} $ is an eigenvalue of multiplicity three; and if $ B_{01}\neq0,\; B_{02}= 0,\;\text{and}\; c_{a}\neq c_{s}  $ then $c_{a}$ as an eigenvalue of multiplicity two. However, if $ B_{01}\neq 0\;\text{and}\;B_{02}\neq 0 $ the wave speeds are distinct and the left side of \eqref{b.2.3} is strictly hyperbolic; here we focus on this case only.

The right eigenvectors of $ \mathbf{A}(\mathbf{U}_0) $ associated with the eigenvalues $\lambda_{i} $ are denoted by $ \mathbf{R}_{i} $, are given by
\begin{equation}\label{b.2.6}
\setlength{\jot}{5pt}
\begin{aligned}
&\mathbf{R}_{1}=(~\rho_{0},~0,~0,~0,
~\left(\frac{\rho P_{\rho}}{P_{s}}\right)_0,~0,~0~)^{tr},\qquad\qquad\qquad\\
&\mathbf{R}_{2,3}=(~0,~0,~0,~\mp c_{a},
~0,~0,~ -B_{01} ~)^{tr},\qquad\qquad\qquad\\
&\mathbf{R}_{4,5}=(~\rho_{0},~\mp c_{s},~\pm \frac{c_{s} B_{01}B_{02}}{\rho_{0}(c_{s}^2-c_{a}^2)},~0,
~0,\frac{B_{02}c_{s}^2}{(c_{s}^2-c_{a}^2)},~ 0 ~)^{tr},\qquad\qquad\qquad\\
&\mathbf{R}_{6,7}=(~\rho_{0},~\mp c_{f},~\pm \frac{c_{f} B_{01}B_{02}}{\rho_{0}(c_{f}^2-c_{a}^2)},~0,
~0,\frac{B_{02}c_{f}^2}{(c_{f}^2-c_{a}^2)},~ 0 ~)^{tr}.\qquad\qquad\qquad\\
\end{aligned}
\end{equation}
The associated left eigenvectors $ \mathbf{L}_{i},\; i = 1, 2,\ldots,7 ,$ can be obtained using the normalization condition
$  \mathbf{L}_{i}. \mathbf{R}_{i} = \delta_{ij}$ where $ \delta_{ij}  $ is Kronecker delta.

We look for a small amplitude high frequency asymptotic solution of the system \eqref{b.2.3} of the form 
\begin{equation}\label{b.2.9}
\mathbf{U}= \mathbf{U}_{0}+
\epsilon^{a_{1}} \mathbf{U}_{1}(x,t,\xi,\tau)+
\epsilon^{a_{2}} \mathbf{U}_{2}(x,t,\xi,\tau)+
\epsilon^{a_{3}} \mathbf{U}_{3}(x,t,\xi,\tau)+
...
\end{equation}
where $\, \displaystyle {\xi=\frac{x}{\epsilon^{b}},~ \tau=\frac{t}{\epsilon^{e}}}\, $ are fast
variables, $ a_{1}<a_{2}<a_{3} $ and $ b,e $ are positive real numbers to be specified later, $ \epsilon $ is a small parameter ($ 0<\epsilon \ll 1$) brought into the problem through initial or boundary condition and is regarded as the strength of the perturbed disturbance. It is the ratio of a typical wavelength relative to the wave modulation length scale and also the ratio of dimensioned wave amplitude relative to a parameter with the same dimension appearing in the problem.  

Each dissipative mechanism existing in the flow defines a local characteristic length (or time)  scale. Short wave or high frequency wave assumption is based on the fact that the wavelength of the wave is much smaller than any other characteristic length scale in the problem. In order to incorporate both dispersive and dissipative effects, we assume that dispersion coefficient $ \chi  $  is larger in magnitude than the diffusion coefficients, such that  
\begin{equation}\label{b.2.7}
\quad\chi = \epsilon^{g}\hat{\chi},\quad \kappa=\epsilon^h \hat{\kappa},\quad
\eta = \epsilon^h \hat{\eta},\quad
\mu= \epsilon^h\hat{\mu},\quad\quad\quad\quad\quad\quad\quad\quad\quad\quad
\end{equation} 
where $g<h$ are positive real numbers to be specified later; the hats designate order one parameters. Our scaling implies that
\begin{equation}\label{b.2.8}
\mathbf{M}_{d}=\epsilon^{g}\,\widehat{\mathbf{M}}_{d},\qquad\mathbf{M}_{v}=\epsilon^{h}\,\widehat{\mathbf{M}}_{v}.\quad\quad\quad\quad\quad\quad\quad\quad\quad\quad\quad\quad\quad\quad\quad\quad
\end{equation}
\noindent
Furthermore, the coefficients matrices $ \mathbf{A},\,\mathbf{M},\,\mathbf{N}\;\text{and}\;\mathbf{Q}\,$ can be expended in a Taylor series about the constant state $ \mathbf{U}_{0} $ as
\begin{equation}\label{b.2.10}
\mathbf{Y}(\mathbf{U})= \mathbf{Y}(\mathbf{U}_{0})+ \nabla\mathbf{Y}(\mathbf{U}_{0}).
(\epsilon^{a_{1}} \mathbf{U}_{1}+\epsilon^{a_{2}} \mathbf{U}_{2}+\ldots)+
\frac{1}{2}\nabla^2\mathbf{Y}(\mathbf{U}_{0}).(\epsilon^{a_{1}} \mathbf{U}_{1}+\ldots)(\epsilon^{a_{1}} \mathbf{U}_{1}+\ldots)^{tr}+\ldots
\end{equation}
where $\mathbf{Y} $ may represent any of the matrices 
$ \mathbf{A},\,\mathbf{M},\,\mathbf{N}\;\text{and}\;\mathbf{Q}\,.$
Now using the derivative transformation
$ \displaystyle{\frac{\partial}{\partial t} = \frac{\partial}{\partial t} + \epsilon^{-e}\frac{\partial}{\partial \tau},\quad \frac{\partial}{\partial x} = \frac{\partial}{\partial x} + \epsilon^{-b}\frac{\partial}{\partial \xi}  }, $ 
\eqref{b.2.8}, \eqref{b.2.10} and perturbation expansion  \eqref{b.2.9}  in \eqref{b.2.3} yield the equation
\begin{equation}\label{b.2.90}
\begin{aligned}
&\biggl[\,\frac{\partial}{\partial t} + \epsilon^{-e}\frac{\partial}{\partial \tau}+
[\,\mathbf{A}_{0}+(\epsilon^{a_{1}} \mathbf{U}_{1}+\epsilon^{a_{2}} \mathbf{U}_{2}+\ldots\,). \nabla\mathbf{A}_{0}+\ldots\,]
\,.\left(\frac{\partial}{\partial x} + \epsilon^{-b}\frac{\partial}{\partial \xi }\right) \\
&
-[\,\epsilon^{g}\,(\widehat{\mathbf{M}}_{0d}+\ldots) + \epsilon^{h}\,(\widehat{\mathbf{M}}_{0v}+\ldots)~]\,.
\left(\frac{\partial^2}{\partial x^2} + 2 \epsilon^{-b}\frac{\partial^2}{\partial \xi \partial x } + 
\epsilon^{-2 b}\frac{\partial^2}{\partial \xi^2 } 
\right)
-[\,\mathbf{N}_{0}+ \ldots\,]\,.
\\
&
\left(\frac{\partial}{\partial x} + \epsilon^{-b}\frac{\partial}{\partial \xi }\right) 
[\,\mathbf{Q}_{0}+ \ldots\,]
\,.\left(\frac{\partial}{\partial x} + \epsilon^{-b}\frac{\partial}{\partial \xi }\right)
\biggr]
(\mathbf{U}_{0}+\epsilon^{a_{1}} \mathbf{U}_{1}+\epsilon^{a_{2}} \mathbf{U}_{2}+\epsilon^{a_{3}} \mathbf{U}_{3}+\ldots)=0.
\end{aligned}
\end{equation} 
To get both dissipation and dispersion effects into picture, we take $ b=e=a_{1}=1,~ g=a_{2}=3/2,~ \text{and }~ h=a_{3}=2$. Equating the coefficients of $ \epsilon^{n/2} $ to zero; this leads to the following system of PDEs satisfied by $ U_1,~  U_2 $, and $ U_3 $
\begin{equation}\label{b.2.11}
\setlength{\jot}{5pt}
\begin{aligned}
&\quad\mathbf{U}_{1\tau} + \mathbf{A}_{0}\mathbf{U}_{1\xi}=0,\qquad\\
&\quad\mathbf{U}_{2\tau} + \mathbf{A}_{0}\mathbf{U}_{2\xi}= \widehat{\mathbf{M}}_{0d}\mathbf{U}_{1\xi\xi},\qquad\\
&\quad\mathbf{U}_{3\tau} + \mathbf{A}_{0}\mathbf{U}_{3\xi} +
\mathbf{U}_{1t} + \mathbf{A}_{0}\mathbf{U}_{1x} +
\nabla\mathbf{A}_{0}.\mathbf{U}_{1\xi}\mathbf{U}_{1\xi} = \widehat{\mathbf{M}}_{0v}\mathbf{U}_{1\xi\xi}+
\widehat{\mathbf{M}}_{0d}\mathbf{U}_{1\xi\xi\xi}.
\qquad
\end{aligned}
\end{equation}
Where $\; \mathbf{A}_{0}=\mathbf{A}(\mathbf{U})_{0},\;
\nabla\mathbf{A}_{0}=\nabla\mathbf{A}(\mathbf{U}_{0}),\;
\widehat{\mathbf{M}}_{0v}=\widehat{\mathbf{M}}_{0v}(\mathbf{U}_{0}),
~\text{and}~ 
\widehat{\mathbf{M}}_{0d}=\widehat{\mathbf{M}}_{0d}(\mathbf{U}_{0}),\;$

The solution of \eqref{b.2.11}(i) is given by 
\begin{equation}\label{b.2.12}
\mathbf{U}_{1}=\sum_{j=1}^{7}\sigma_{j}(x,t,\theta_{j})\,\mathbf{R}_{j},
\end{equation}
where $ \sigma_{j} = (\mathbf{L}_{j}.\mathbf{U}_{1}) $ is an arbitrary scalar valued function called the wave amplitude; it depends on the $ j $-th phase variable $ \theta_{j} $ given by
$ \theta_{j} = k_{j}\xi -\omega_{j}\tau,$ where the wavenumber $ k_{j} $ and frequency 
$\omega_{j}$ satisfy $\; \omega_{j}=\lambda_{j}k_{j},\; j=1,\ldots,7; $ indeed, the dependence of $ \sigma_{j} $ on $ \theta_{j} $ describes the waveform.

We also assume that $ \sigma_{j}(x,t,\theta) $ has zero mean with respect to the phase variable $ \theta_{j} ,$ i.e., 
\begin{equation*}
\lim_{T\to\infty}\frac{1}{T}\int_{0}^{T}\sigma_{j}(x,t,\theta)~d\theta=0;
\end{equation*}
and frequency-wave number pairs are related by
\begin{equation}\label{b.2.13}
\omega_{j}=\mu_{jnm}\omega_{m}+ \mu_{jmn} \omega_{n},\qquad
k_{j}=\mu_{jnm}k_{m}+\mu_{jmn}k_{n},\qquad\qquad\qquad
\end{equation}
where $$\displaystyle{ \mu_{jmn} =  \frac{k_{j}(\lambda_{j}-\lambda_{m})}{k_{n}(\lambda_{n}-\lambda_{m})}}.\qquad\qquad\qquad\qquad\qquad\qquad\qquad$$
Using \eqref{b.2.12} in \eqref{b.2.11}(ii) and solving the resulting equation for $\mathbf{U}_{2}$ by the method of characteristics, we get
\begin{equation}\label{b.2.14}
\mathbf{U}_{2}=\sum_{j=1}^{7}\sigma_{j\theta}(x,t,\theta_{j})\,\mathbf{P}_{j},
\end{equation}
where the vector $\mathbf{P}_{j}  $ satisfies 
\begin{equation}\label{b.2.15}
(k_{j}\mathbf{A}_{0}-\omega_{j}\mathbf{I})\,\mathbf{P}_{j}=k_{j}^2\,\widehat{\mathbf{M}}_{0d}\,\mathbf{R}_{j},
\end{equation}
where the dispersion matrix  $ \widehat{\mathbf{M}}_{0d} $ satisfies the solvability condition $ \mathbf{L}_{j}.\widehat{\mathbf{M}}_{0d}\,\mathbf{R}_{j}=0.$ Now we use \eqref{b.2.14},  \eqref{b.2.12} in \eqref{b.2.11}(iii) and solve for $ \mathbf{U}_{3}; $
on using the secularity condition  
\begin{equation}\label{b.2.16}
\lim_{\tau\to\infty}\frac{1}{\tau}\mathbf{U}_{3}(\xi,\tau)=0,
\end{equation}
we arrive at the following set of 
following integro-differential equations for the wave amplitudes
$ \sigma_{1}(x,t,\theta),\ldots, \sigma_{7}(x,t,\theta) $ 
\begin{equation}\label{b.2.17}
\begin{aligned}
&\sigma_{jt}(x,t,\theta) +\lambda_{j}\,\sigma_{jx}(x,t,\theta) +
E_{j}\,\sigma_{j}(x,t,\theta)\, \sigma_{j\theta}(x,t,\theta)\\[2pt]
&+\sum_{m < n}^{(j)}\mu_{nmj}\Gamma_{jmn}\lim_{T\to\infty}\frac{1}{T}
\int_{0}^{T}\sigma_{m}(x,t,\mu_{mnj}\theta+\mu_{mjn}\zeta)
\sigma_{n\zeta}(x,t,\zeta)\,d\zeta=
\Omega_{j}\sigma_{j\theta\theta}+\Lambda_{j}\sigma_{j\theta\theta\theta},
\end{aligned}
\end{equation}
where the coefficients are 
\begin{equation}\label{b.2.18}
\setlength{\jot}{5pt}
\begin{aligned}
&E_{j}=k_{j}\,\mathbf{L}_{j}\,.\,\nabla\mathbf{A}_{0}\,.\,\mathbf{R}_{j}\,\mathbf{R}_{j},\\
&\Gamma_{jmn}=\mu_{jmn}\,k_{n}\,\mathbf{L}_{j}\,.\,\nabla\mathbf{A}_{0}\,.\,\mathbf{R}_{m}\,\mathbf{R}_{n}+
\mu_{jnm}\,k_{m}\,\mathbf{L}_{j}\,.\,\nabla\mathbf{A}_{0}\,.\,\mathbf{R}_{n}\,\mathbf{R}_{m},\\
&\Omega_{j} = k_{j}^2\,\mathbf{L}_{j}\,.\,\widehat{\mathbf{M}}_{0v}\,\mathbf{R}_{j},\qquad
\Lambda_{j} = k_{j}^2\,\mathbf{L}_{j}\,.\,\widehat{\mathbf{M}}_{0d}\,\mathbf{P}_{j}.
\end{aligned}
\end{equation}
Here, we have seven integro-differential equations with dissipative and dispersive terms corresponding to the MHD waves in the real gas background. In each case, except for the entropy wave, we have a non zero (for the real gas we considered) self interaction coefficient $ E_{j} $ called the nonlinearity parameter. For the entropy wave, it is zero which can be attributed to the linearly degenerate behavior of the entropy wave \cite{MR760229}; also the entropy wave equation has no integral term, since it is not influenced by the interaction of the other waves due to the fact that entropy wave is a Riemann invariant. The linear integral term in  other cases results from the three wave interactions between different waves and its coefficient corresponds to the amount of wave produced through the interaction between other two waves. 
\section{Resonant interaction between fast magnetoacoustic entropy triad}\label{evo}
In this section, the general evolution equation \eqref{b.2.17} derived in the last section is restricted to the study of three wave interaction at a time. The resulting triads can be classified as $ (a) $  magnetoacoustic and  entropy waves,  $ (b) $  Alfvén and entropy waves, $ (c) $ Alfvén and magnetoacoustic waves, and $ (d) $ slow and fast magnetoacoustic waves.

We focus our study on the two fast magnetoacoustic and entropy wave interactions and obtain the precise formulation of the transport equation in this case from the general 
equation \eqref{b.2.17}. For this purpose, we make certain assumptions: $ (i) $ the wave amplitudes are $ 2\pi$-periodic function of the phase $ \theta ,$ 
$ (ii) $ all other waves except fast magnetoacoustic and entropy wave triad are in non-resonance (hence the integral term corresponding to them become zero), and $ (iii) $
the fundamental harmonics of all waves in our triad satisfy the resonance conditions \begin{equation}\label{b.2.19}
\omega_{1}+\omega_{6}+\omega_{7}=0,\qquad\qquad\qquad k_{1}+k_{6}+k_{7}=0;
\end{equation} 
which are satisfied when $ k_{6}=k_{7}\; \text{and}\; k_{1}=-2k_{7},$ where the wave numbers $k_{1},\,k_{6},\,k_{7} $ correspond to the entropy wave and left and right moving magnetoacoustic waves, respectively. Using the restrictions $ (i),(ii)~\text{and}~ (iii),$ we get the following system of asymptotic equations

\begin{eqnarray}\label{b.2.20}
\begin{aligned}
&\sigma_{6t}- c_{f}\sigma_{6x} -
k_{6}E_{f}\sigma_{6}\sigma_{6\theta} -
k_{6}\mu_{761}M_{f}\frac{1}{2\pi}
\int_{0}^{2\pi}\sigma_{1}(-\theta-\zeta)
\sigma_{6\zeta}(\zeta)\,d\zeta=
k_{6}^2\Omega_{f}\sigma_{6\theta\theta}-k_{6}^3\Lambda_{f}\sigma_{6\theta\theta\theta},\\
&\sigma_{1t}(\theta)=k_{1}^2\Omega_{e}\sigma_{1\theta\theta},\\
&\sigma_{7t} + c_{f}\sigma_{7x} +
k_{7}E_{f}\sigma_{7}\sigma_{7\theta} +
k_{7}\mu_{167}M_{f}\frac{1}{2\pi}
\int_{0}^{2\pi}\sigma_{6}(-\theta-\zeta)
\sigma_{1\zeta}(\zeta)\,d\zeta=
k_{7}^2\Omega_{f}\sigma_{7\theta\theta}+k_{7}^3\Lambda_{f}\sigma_{7\theta\theta\theta},\\
\end{aligned}
\end{eqnarray}
where the coefficients are 
\begin{equation}\label{b.2.21}
\setlength{\jot}{5pt}
\begin{aligned}
&G_{0}=\left(1+\frac{\rho_{0}}{c_{0}}c_{\rho 0}\right),\qquad 
H_{0}=\frac{1}{2}\frac{\rho_{0}P_{\rho s 0}}{P_{s0}},\qquad 
\gamma_{f}=\frac{c_{f}^2-c_{0}^2}{c_{f}^2-c_{s}^2}, \qquad 
\gamma_{s}=\frac{c_{s}^2-c_{0}^2}{c_{s}^2-c_{f}^2},\\
& E_{f}=\left(G_{0}\,\gamma_{s} + \frac{3}{2} \,\gamma_{f}\right)c_{f}, 
\qquad M_{f}=\left((G_{0}-H_{0})\,\gamma_{s} + \frac{\gamma_{f}}{2}\right)c_{f} \,, \\
&\Omega_{f} =\frac{c_{a}^2}{2\rho_{0}}\left(\frac{4\,\gamma_{s}}{3\,c_{s}^2}  +
\frac{\gamma_{f} }{c_{f}^2} \right)\hat{\mu}
+\left(\frac{\gamma_{f}\,}{2}\right)\hat{\eta}+
\left(\frac{\gamma_{s}\delta}{2\rho_{0}c_{p0}}\right)\hat{\kappa},
\qquad
\Omega_{e} =\left(\frac{1}{\rho_{0}c_{p0}} \right)\hat{\kappa},\\
&\Lambda_{f} = \left(\frac{c_{a}^2}{c_{f}^2-c_{a}^2}\right)\rho_{0}\gamma_{f}c_{f}\hat{\chi}^2.
\end{aligned}
\end{equation} 
These equations, with the viscous and dispersive terms are similar in form to the equations for wave interaction in elasticity \cite{MR1160152} and gas dynamics \cite{MR760229}. However, the MHD waves have the dispersive behavior which is not present in acoustic and elastic waves and they encapsulate the KdV type behavior in addition to the resonant interaction. The dynamics of combination of resonant interaction and weak dispersion is complicated. To study the behaviour of these two effects we consider that there are no spacial modulation and we also consider that the viscosity is absent;
consequently the second equation implies that the entropy is independent of time and is given by the initial condition. Hence, our system is reduced to a pair for fast magnetosonic waves and entropy wave interaction.  Introducing the new independent variables, 
\begin{equation}\label{b.2.22}
\setlength{\jot}{5pt}
\begin{aligned}
&\alpha_{1}(x,t,\theta)=k_{6}\sigma_{6}(x,t,\theta),\\
&\alpha_{2}(x,t,\theta)=k_{7}\sigma_{7}(x,t,\theta),\\
& K(x,\theta)= k_{1}\sigma_{1\theta}(x,t,\theta),
\end{aligned}
\end{equation}
the resulting pair of equations  is given by
\begin{equation}\label{b.2.23}
\begin{aligned}
&\alpha_{1t} +
E_{f}\alpha_{1}\alpha_{1x} +
M_{f}\frac{1}{2\pi}
\int_{0}^{2\pi}K(\zeta-x)\,
\alpha_{2}(\zeta)\,d\zeta+
\Lambda_{f}\alpha_{1xxx}=0,\\
&\alpha_{2t} +
E_{f}\alpha_{2}\alpha_{2x} -
M_{f}\frac{1}{2\pi}
\int_{0}^{2\pi}K(x-\zeta)\,
\alpha_{1}(\zeta)\,d\zeta+
\Lambda_{f}\alpha_{2 xxx}=0,\\
\end{aligned}
\end{equation}
where a change of variables from $ \theta $ to $ x $ and $ \alpha_{1}(x,t) $ to $ \alpha_{1}(-x,t) $
is used. Finally, if we take the amplitudes $ \alpha_{1}=\alpha_{2} $ and odd kernel $ K(-x)=-K(x) $
our system \eqref{b.2.23} is symmetric with respect to $ \alpha_{1} $ and $ \alpha_{2}$ and it is reduces to the following single equation
\begin{equation}\label{b.2.24}
\alpha_{t} +
E_{f}\,\alpha \,\alpha_{x} +
M_{f}\,\frac{1}{2\pi}
\int_{0}^{2\pi}K(y-x)\,
\alpha(y)\,dy+
\Lambda_{f}\,\alpha_{xxx}=0.
\end{equation}

Now, we introduce dimensionless variables, defined as;
\begin{equation*}
\begin{aligned}
&x^{*}=\frac{x}{l},\;\; t^{*}=\frac{t\sqrt{P_{0}/ \rho_{0}}}{l},
\;\; P^{*}=\frac{P}{P_{0}}, \;\; \rho^{*}=\frac{\rho}{\rho_{0}},
\;\; T^{*}=\frac{T}{T_{0}}, \;\; s^{*}=\frac{s}{s_{0}},
\;\; b^{*} = b \rho_{0},\\
&\;\; u^{*}=\frac{u}{\sqrt{P_{0}/ \rho_{0}}},
\;\; v^{*}=\frac{v}{\sqrt{P_{0}/ \rho_{0}}},
\;\; w^{*}=\frac{u}{\sqrt{P_{0}/ \rho_{0}}},
\;\; B_{1}^{*}=\frac{B_{1}}{B_{0}},
\;\; B_{2}^{*}=\frac{B_{2}}{B_{0}},
\;\; B_{3}^{*}=\frac{B_{3}}{B_{0}};
\end{aligned}
\end{equation*}
using these variables in \eqref{b.2.24}, and dropping the $ * $ sign, the equation remain unchanged in the dimensionless form. Using equation of state \eqref{b.2.2}, the coefficients evaluated at the undisturbed state are given by
\begin{equation}\label{b.2.25}
\setlength{\jot}{5pt}
\begin{aligned}
&c_{0}=\sqrt{\frac{(\delta+1)}{(1-b)}},\qquad 
G_{0}=\frac{(\delta+2)}{2(1-b)},\qquad 
H_{0}=\frac{(\delta+1)}{2(1-b)},\quad 
\gamma_{f}=\frac{c_{f}^2-c_{0}^2}{c_{f}^2-c_{s}^2}, \quad 
\gamma_{s}=\frac{c_{s}^2-c_{0}^2}{c_{s}^2-c_{f}^2},\\
& E_{f}=\left(G_{0}\,\gamma_{s} + \frac{3}{2}\right)c_{f} \,\gamma_{f}, 
\qquad M_{f}=\left((G_{0}-H_{0})\,\gamma_{s} + \frac{\gamma_{f}}{2}\right)c_{f} \,,\quad 
\Lambda_{f} = \left(\frac{c_{a}^2}{c_{f}^2-c_{a}^2}\right)\gamma_{f}c_{f}\hat{\chi}^2.\\\notag
\end{aligned}
\end{equation} 

\section{Numerical solutions}\label{num1}
In this section, we study the numerical solutions of the transport equation \eqref{b.2.24} and discuss the effects of various parameters in the light of the real gas background. Zabusky and Kruskal \cite{MR716190} have studied the solitary wave solution for the KdV-type equation
\begin{equation}\label{key}
u_{t}+uu_{x}+\beta^2u_{xxx}=0, 
\end{equation} with a small amount of dispersion $ \beta=0.022 $ and periodic initial data. We investigate the existence of solitary wave solutions of \eqref{b.2.24} in the light of real gas and magnetic field effects. For numerical computation, we use pseudo-spectral method developed by Fronberg and Whitham \cite{MR497916}, which is suitable for the evaluation of certain operators, and can considerably speed up the calculations while using fast Fourier transform. We use the trapezoidal rule to evaluate the integral and the temporal variable is discretized using a fourth-order Runge-Kutta method. The evolution equation \eqref{b.2.24} is composed of a Burgers equation type nonlinear term with a small amount of dispersion coefficient ($ \Lambda_{f}\approx 0.0013 $ ) and an integral term which is weakly dispersive \cite{MR1719749}. 

In numerical experiments, we have taken the convolution with the kernel $ K(x)= \sin x $ and the periodic initial data $ \alpha(x,0) = A \cos x ,$ with $ A=1,2. $
Since the dispersive effects are negligible, the corresponding terms can be neglected in
\eqref{b.2.24} leading to $ \alpha_{t} + E_{f}\, \alpha\, \alpha_{x}=0 $; 
it's solution is given by the implicit equation $ \alpha =  A \cos (x- E_{f} \alpha t) $. 

Initially the nonlinear term dominates the solution, and a steepening is seen in the region of negative slope, we find that $ \alpha $ tends to become discontinuous at the breakdown time $ t_{b}=\left( \frac{1}{A\, E_{f}}\right) $ but after some time the dispersive term becomes dominant and instead of discontinuity, oscillations of small wavelength develop on the left part of the wave profile;  their amplitudes increase, and each oscillation achieves steady amplitude after some time. We consider the effects of van der Waals parameter $ b $ and magnetic field  on the formation of these solitons.

\begin{figure}[!t]
	\centering
	\includegraphics[width=1.00\textwidth]{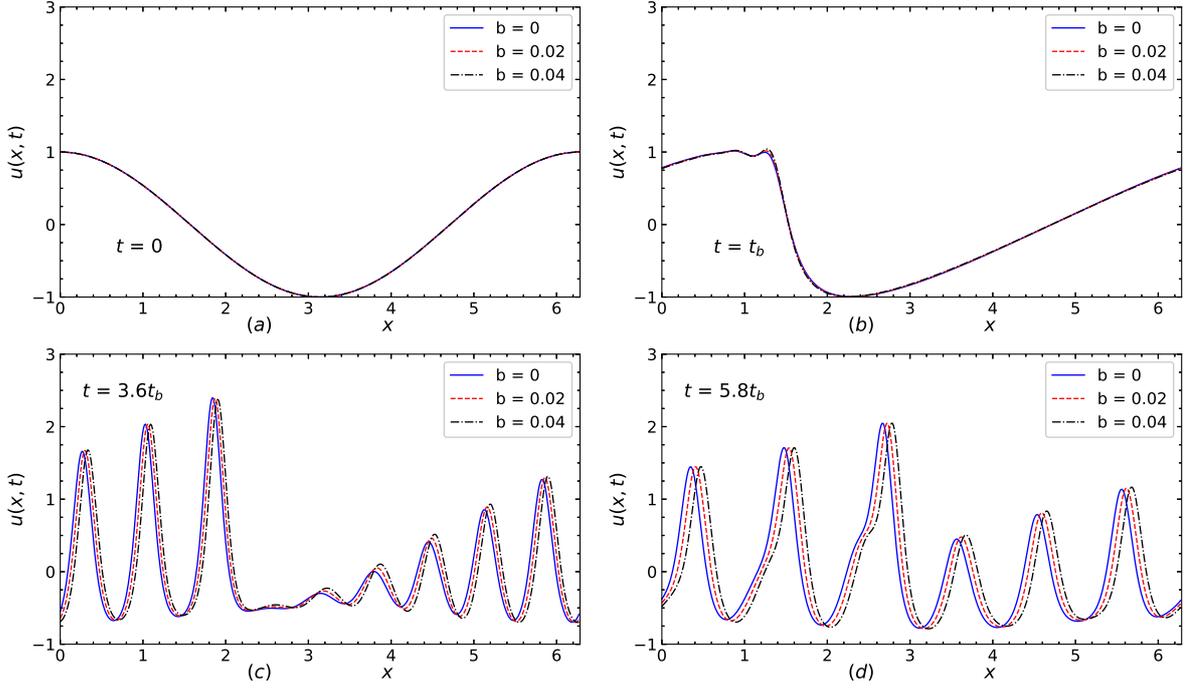}
	\caption{Case(I)- Three families of solitary wave profiles for $b=0,\,b=0.02,\,~\text{and}~ b=0.04 $  with initial data $ a(x,0)= \cos x $ at different times, $ t_{b} $ is the breakdown time,  $ \delta =0.04 ,\chi=1,\, B_{01}= 0.1,\,\text{and}\, B_{02}=1.$  }
	\label{fig41}
\end{figure}

\begin{figure}[!t]
	\centering
	\includegraphics[width=1.00 \textwidth]{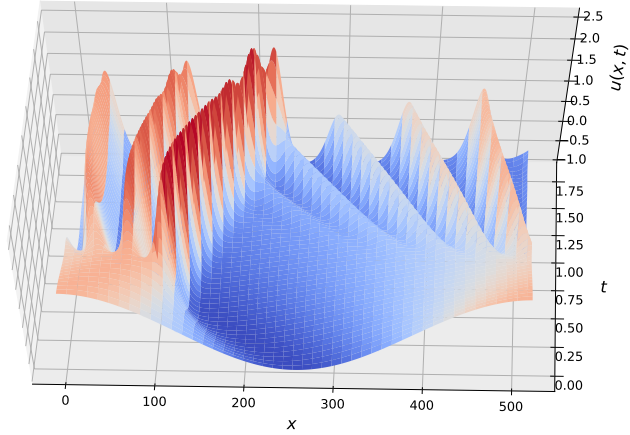}
	\caption{Space-time (top view) of the temporal development of solitons  with initial data $ a(x,0)= \cos x $ for  $ \delta =0.04 ,\,\chi=1,\, B_{1}= 0.1,\,\text{and}\, B_{02}=0.1$. }
	\label{fig42}
\end{figure}

For case-(I), we perform three sets of numerical experiments each for van der Waal parameter
$ b=0,\,b=0.02,\,\text{and}\, b=0.04 $ with initial data $ \alpha (x,0)=\cos x $ and $ \delta =0.04 ,\,\chi=1,\, B_{01}= 0.1,\, B_{02}=0.1.$ Fig. \ref{fig41}(b) shows the wave profiles at breaking time $ t_{b} $, which is almost same in all cases as depicted in Table (\ref{t1.1}); there appears a minute oscillatory behavior in each case which is due to neglected dispersive term. Solitary wave train formations and their overlapping can be seen in  Fig. \ref{fig41}(c) and Fig. \ref{fig41}(d), respectively; a space-time evolution of the wave profiles is displayed in (\ref{fig42}).

\begin{table}[b]
	\caption {Evolution equation coefficients for case-(I).}\label{t1.1} 
	\begin{center}
		\begin{tabular}{ |c|c|c|c|c|c|c| } 
			\hline
			& $ E_{f} $ & $ M_{f}$ & $ \Lambda_{f}$ & $ t_{b}$ \\ 
			\hline
			b=0 & 2.7756 & 1.0134 & 0.00139 & 0.360 \\ 
			b=0.02&2.8096 & 1.0278 & 0.00137 & 0.355 \\ 
			b=0.04& 2.8457 & 1.0430 & 0.00136 & 0.351 \\ 
			\hline
		\end{tabular}
	\end{center}
\end{table}

\begin{figure}[!t]
	\centering
	\includegraphics[width=1.00\textwidth]{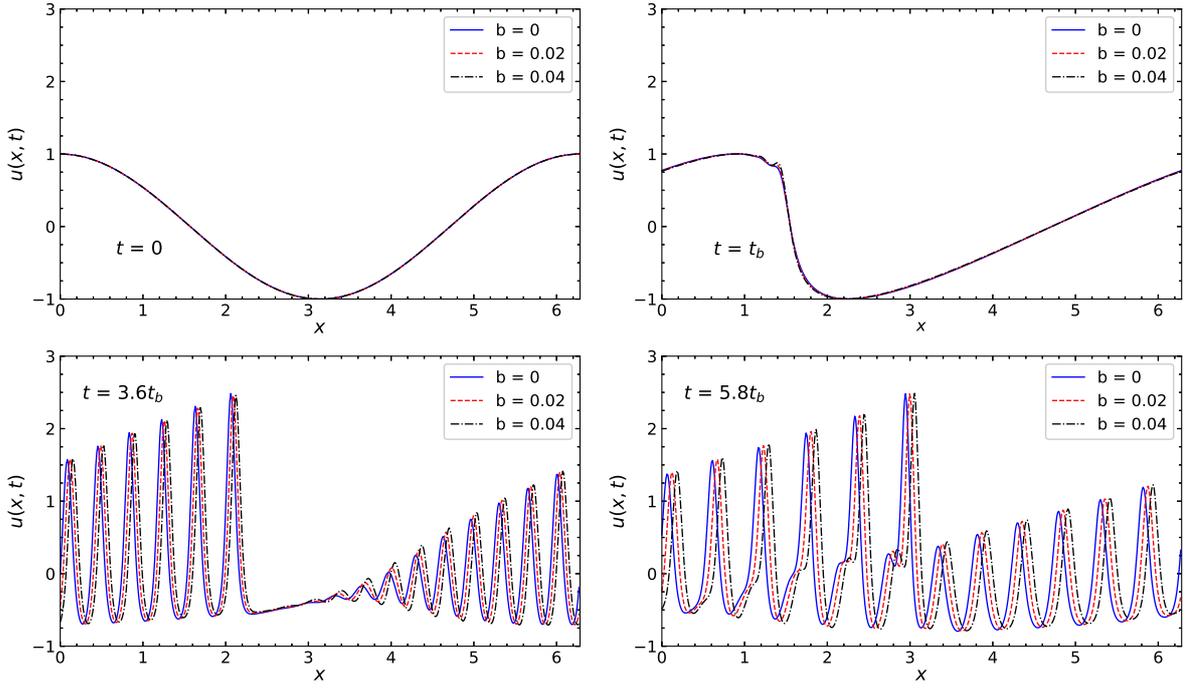}
	\caption{Case(II)-Three families of solitary wave profiles for $b=0,\,b=0.02,\,~\text{and}~ b=0.04 $ with initial data $ a(x,0)= \cos x $ at different times, $ t_{b} $ is the breakdown time,  $ \delta =0.04 ,\,\chi=1,\, B_{01}= 0.05,\,\text{and}\, B_{02}=1$}
	\label{fig43}
\end{figure}

For case-(II), we perform three sets of experiments with the same data as in case-(I) except the magnetic field parameter, which is $ B_{01}=0.05$. In this case, we find that the breakdown time is almost same as in the previous case depicted in Fig.(\ref{fig43})(b) and Table (\ref{t1.2}). Fig.(\ref{fig43})(c) shows the soliton formation; however the number of solitary waves increases and their width decreases as compared to the previous case, which can be explained by the fact that there is considerable decrease in the dispersion coefficient $ \Lambda  $ with the decrease in $ B_{01} $ as shown in Table (\ref{t1.2}); Zabusky and Kruskal \cite{MR716190} had shown that, thickness of the pulse decreases with the decrease in the dispersion parameter and increase in the wave amplitude.
\begin{table}[b]
	\caption {Evolution equation coefficients for case-(II).}\label{t1.2} 
	\begin{center}
		\begin{tabular}{ |c|c|c|c|c|c|c| } 
			\hline
			& $ E_{f} $ & $ M_{f}$ & $ \Lambda_{f}$ & $ t_{b}$ \\ 
			\hline
			b=0 & 2.7725 & 1.0126 & 0.00034 & 0.366 \\ 
			b=0.02&2.8066 & 1.0271 & 0.00034 & 0.356 \\ 
			b=0.04& 2.8428 & 1.0423 & 0.00033 & 0.351 \\ 
			\hline
		\end{tabular}
	\end{center}
\end{table}

\begin{figure}[!t]
	\centering
	\includegraphics[width=1.00\textwidth]{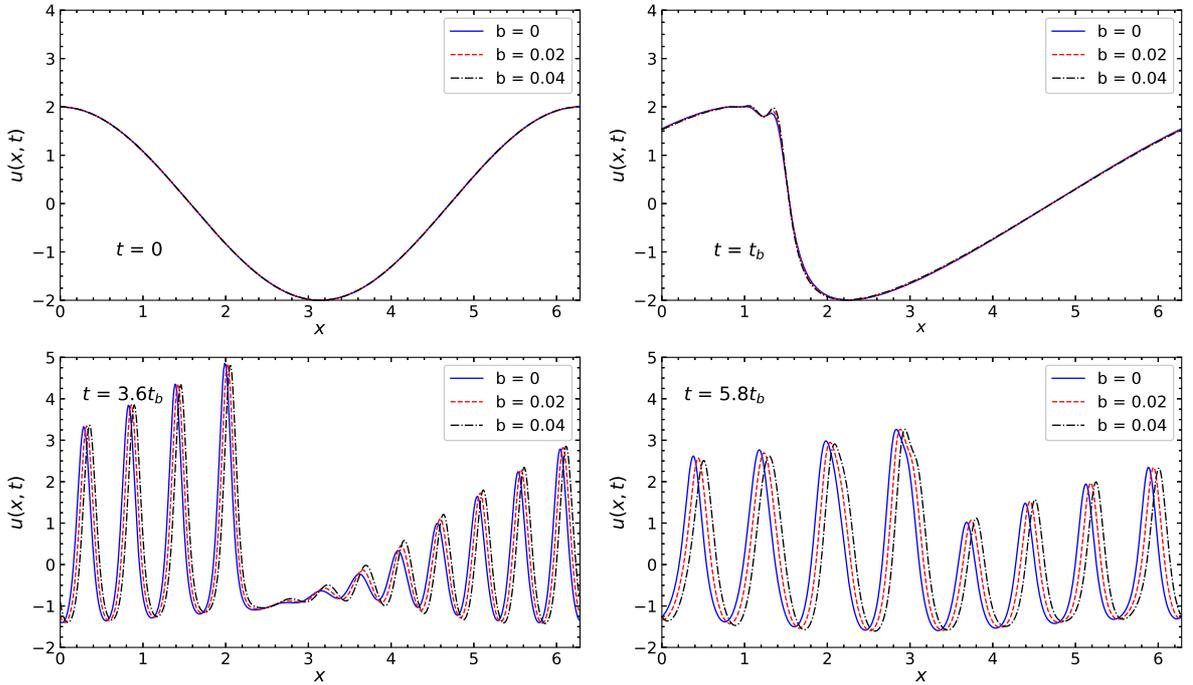}
	\caption{Case(III)-Three families of solitary wave profiles for $b=0,\,b=0.02,\,~\text{and}~ b=0.04 $ with initial data $ a(x,0)= 2 \cos x $ at different times, $ t_{b} $ is the breakdown time,  $ \delta =0.04 ,\,\chi=1,\, B_{01}= 0.1,\,\text{and}\, B_{02}= 1$ }
	\label{fig44}
\end{figure}

Numerical solution for case-(III), with the same set of parameters as in case-(I) except with initial condition $ a(x,0)=2 \cos x $, i.e., with twice the amplitude, are obtained with three sets of experiments. In this case, as is clear from the expression of $ t_{b} $ that the breakdown time decreases to one half that of the first case. Indeed, as the amplitude of the initial profile increases, the amplitude of the solitary wave increases, while, the number of pulses increases and their thickness decreases as depicted in the Figure (\ref{fig44}), which is again explained by the result mentioned in the last case.

\begin{figure}[!t]
	\centering
	\includegraphics[width=1.0\textwidth]{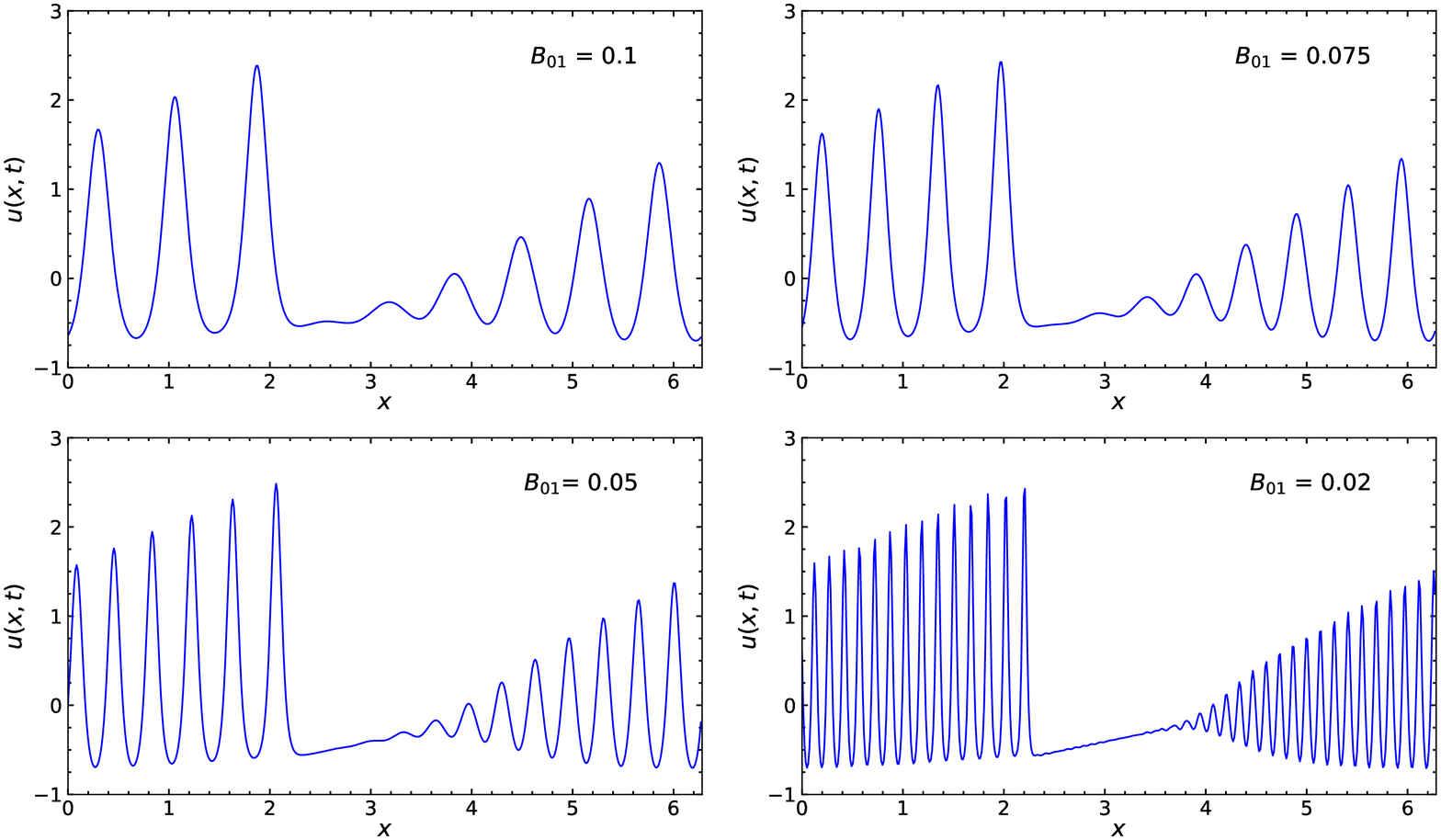}
	\caption{Case(IV)-The solitary wave profiles  with different values of magnetic fields  $ B_{01}= 0.1,\, B_{01}= 0.075,\, B_{01}= 0.05,\,\text{and}\, B_{01}= 0.02, $ with initial data $ a(x,0)= \cos x $ for $b=0.02$,  $ \delta =0.04 ,\,\chi=1,\,\text{and}\, B_{02}=0.1.$ }
	\label{fig45}
\end{figure}

For case-(IV), we investigate the influence of magnetic field $ B_{01} $ on the formation of solitons at same time, with the same set of parameters as in case-(I) but with the decreasing magnetic field $ B_{01} $ which is exhibited in the Figure (\ref{fig45}). The number of solitons increases and their width decreases, which is illustrated by a small value of $ \Lambda_{f} $ shown in Table (\ref{t1.3}); and eventually with the decrease in the magnetic field, there is a breakdown of the solution,  whereas the breakdown time $ t_{b} $ is almost same in each case as displayed in Table (\ref{t1.3}) 

\begin{table}[b]
	\caption {Evolution equation coefficients for case-(IV).}\label{t1.3} 
	\begin{center}
		\begin{tabular}{ |c|c|c|c|c|c|c| } 
			\hline
			& $ E_{f} $ & $ M_{f}$ & $ \Lambda_{f}$ & $ t_{b}$ \\ 
			\hline
			$ B_{01}=0.1 $ & 1.9050 & 0.7465 & 0.00087 & 0.524 \\ 
			$ B_{01}=0.075 $&1.9046& 0.7465 & 0.00049 & 0.524 \\ 
			$ B_{01}=0.05 $& 1.9044 & 0.7465 & 0.00021 & 0.525 \\ 
			$ B_{01}=0.02 $& 1.9042 & 0.7465 & 0.00011 & 0.525 \\ 
			\hline
		\end{tabular}
	\end{center}
\end{table}
\begin{figure}[!t]
	\centering
	\includegraphics[width=1.0\textwidth]{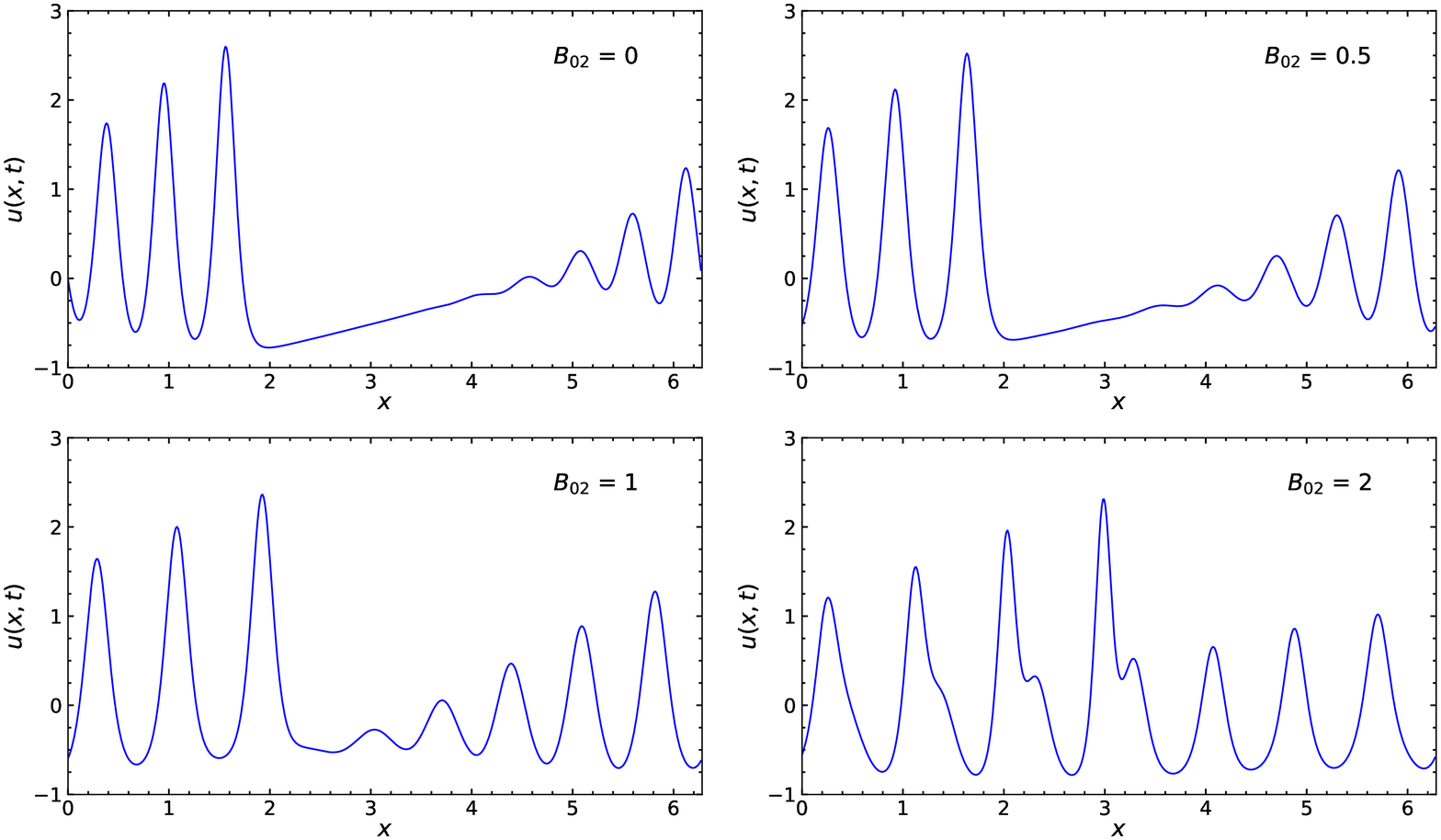}
	\caption{Case(V)-The solitary wave profiles with different values of magnetic fields  $ B_{02}= 0,\, B_{02}= 0.5,\, B_{02}= 1,\,\text{and}\, B_{02}=2, $ with initial data $ a(x,0)= \cos x $ and $b=0.02$,  $ \delta =0.04 ,\,\chi=1,\, \,\text{and}\, B_{01}=0.1.$ }
	\label{fig47}
\end{figure}
Finally, in case-(V) and case-(VI), we consider the influence of magnetic field $ B_{02} $ on the solitary wave solutions with the same set of parameters as in case-(I) but with the increasing magnetic field. For case-(V) in Figure (\ref{fig47}), we display the solitary wave profiles with increasing magnetic field  and it is noticed  that with the increase in the magnetic field there is a higher development in the solitory wave profile, which is explained by the decrease in the breaking time $ t_{b} $ with the increase in magnetic field as depicted in Table(\ref{t1.4}).
\begin{table}[b]
	\caption {Evolution equation coefficients for case-(V) and case-(VI).}\label{t1.4} 
	\begin{center}
		\begin{tabular}{ |c|c|c|c|c|c|c| } 
			\hline
			& $ E_{f} $ & $ M_{f}$ & $ \Lambda_{f}$ & $ t_{b}$ \\ 
			\hline
			$ B_{02}=0 $ & 1.8940 & 0.7431 & 0.00086 & 0.527  \\ 
			$ B_{02}=0.5 $&2.1549& 0.8238 & 0.00112 & 0.464   \\ 
			$ B_{02}=1.0 $& 2.8096 & 1.0278 & 0.00137 & 0.355 \\ 
			$ B_{02}=2.0 $& 4.6079 & 1.5999 & 0.00190 & 0.217 \\ 
			\hline
		\end{tabular}
	\end{center}
\end{table}
\begin{figure}[!t]
	\centering
	\includegraphics[width=1.0\textwidth]{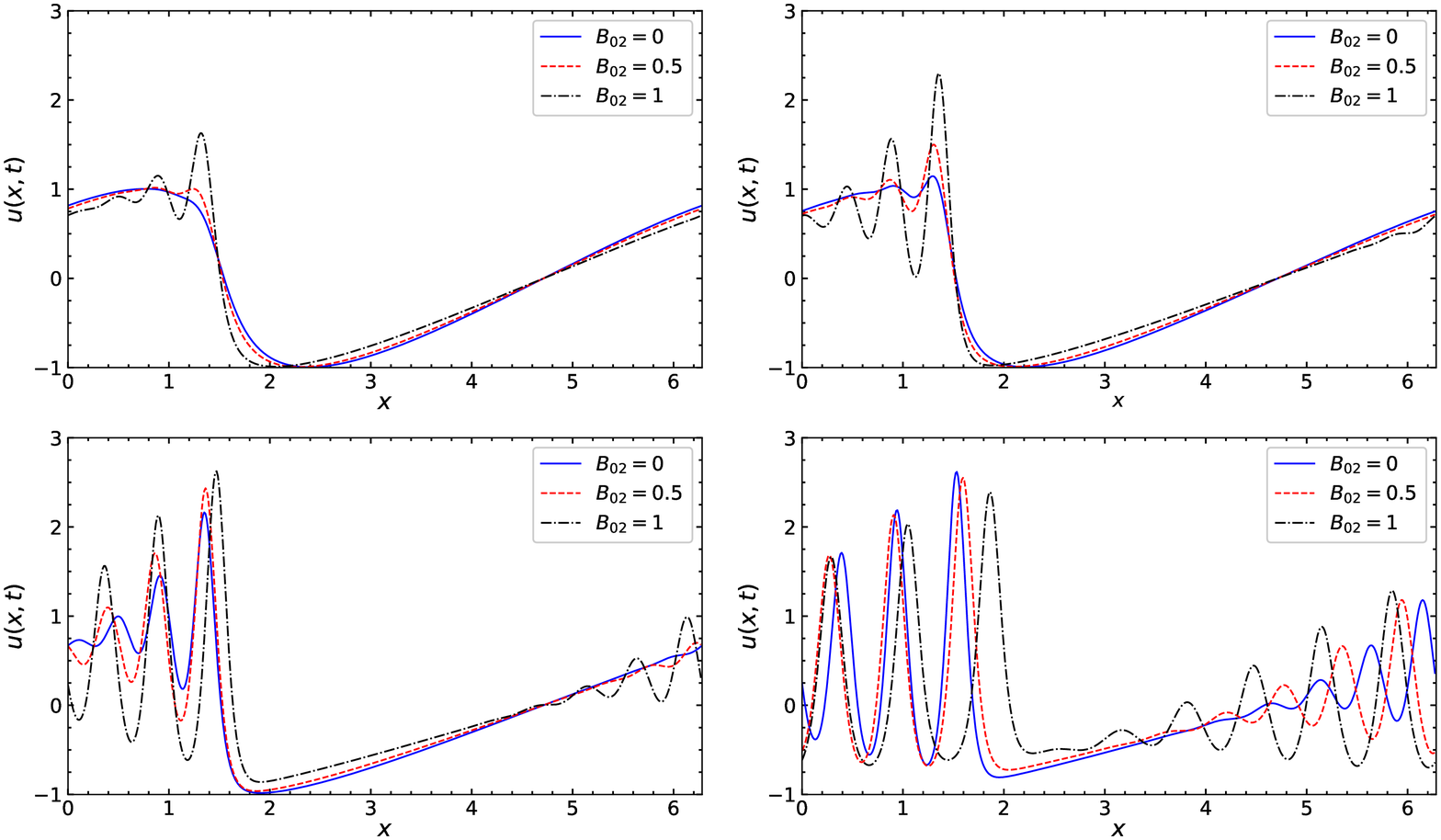}
	\caption{Case(VI)- Three families of solitary wave profiles at the same time for different values of magnetic fields  $ B_{02}= 0,\, B_{02}= 0.5,\, B_{02}= 1,\,\text{and}\,B_{02}= 2, $ respectively, with initial data $ a(x,0)= \cos x $ and $b=0.02$,  $ \delta =0.04 ,\,\chi=1,\, \text{and}\, B_{01}=0.1.$ }
	\label{fig46}
\end{figure}

Similarly, in case-(VI), as portrayed in Figure (\ref{fig46}), three families of wave profiles are displayed at a particular time for three different values of magnetic field; in the first Figure it observed that for a strong magnetic field, i.e., $ B_{02}=1, $ the oscillation begin to set in, whereas for the other two values they are still not seen; with the increase in time the evolution of solitons is always ahead of the wave profiles
that corresponds to the largest amount of magnetic field and this behaviour is due to the smallest value of $ t_{b} $ corresponding to the most significant value of $ B_{02} $ as depicted in Table (\ref{t1.4}).

\newpage
\section{Conclusions}\label{con1}
We have studied, using a perturbation method, a magnetohydrodynamics model in the real gas background with viscosity, magnetic diffusion, thermal conductivity, and magnetic dispersion; we derive a system of evolution equations for the resonant interaction among the characteristics  of MHD. To illustrate the effects of the interplay between the nonlinearity, dispersion, and resonance, we focus on a single triad composed of two opposite moving fast magnetoacoustic waves and an entropy wave with certain assumptions. The resulting single equation has a Burgers type nonlinear term with coefficient $ E_{f} $ (always positive for the real gas considered), dispersion term with the coefficient $ \Lambda_{f} $ and linear integral term which corresponds to the resonant interaction among the waves in considered triad with coefficient $ M_{f} $.
In all the numerical experiments we perform, the dispersion coefficient is small, hence at the initial stage the quadratic nonlinear term dominates the behavior of solution profile, but after some time the dispersion and integral terms come into play and solitary wave formation takes place. 

We investigate, how the real gas and magnetic field effects influence the wave formation, shape and amplitude of the solitary wave profiles for periodic initial data. The effects of van der Waals parameter $ b $ are considered in cases  (I),(II), and (III); solitary wave formation takes places in each case and it is displayed in the Figs. (\ref{fig41}),(\ref{fig43}), and (\ref{fig44}) showing thereby that they have small effect on the breakdown time $ t_{b} $ and on the development of the solitary wave profiles, however an enhancement in the amplitude of the initial wave profile causes the resulting amplitude of the solitary wave profile to increase and to decreases the breakdown time as shown in Fig. (\ref{fig43}). The effects of magnetic field on the wave evolution are considered in cases (II),(IV),(V), and (VI). An increase in the number of soliton and decrease in their thickness with a reduction in the value of $ B_{01} $ is shown in the Fig. (\ref{fig45}) which is explained from the entries of $ \Lambda_{f} $ in the Table (\ref{t1.4}) for different values of $ B_{01} .$ There is a lagging in  evolution of the wave profiles for smaller values of $ B_{02} $ as illustrated in the Fig.(\ref{fig47}); evolution of three different wave profiles (at the same time) corresponding to three different values of $ B_{02} $ is depicted in Fig. (\ref{fig46}) and it is noticed that the breaking time decreases with an increase in the value of $ B_{02} ,$ whereas, the number of solitons and their width remain the same in each case.



\end{document}